 \documentclass[11pt]{article}
\newcommand{\comm}[1]{}

\def\Lindex{\comm}
\def\indexNN{\comm}
\usepackage[pagewise]{lineno}
\usepackage{amssymb}
\usepackage[english]{babel}\usepackage[T1]{fontenc}
\usepackage{lmodern}
\usepackage[T1]{fontenc}
\usepackage[utf8]{inputenc}
\usepackage{amsmath}
\usepackage{amsfonts}
\usepackage{mathtools}
\usepackage{color}

   \csname @addtoreset\endcsname{equation}{section}
\newtheorem{theorem}{Theorem}[section]
\newtheorem{lemma}{Lemma}[section]

\newtheorem{proposition}{Proposition}[section]
\newtheorem{corollary}{Corollary}[section]
\newtheorem{condition}{Condition}[section]

\newtheorem{remark}{Remark}[section]

\def\e{\varepsilon}

\def\defi{\stackrel{{\scriptscriptstyle \Delta}}{=}}

\def\N{\mu}
\def\a{\alpha}

\def\o{\omega}
\def\O{\Omega}

\def\F{{\cal F}}
\def\w{\widehat}
\def\ww{\widetilde}
\def\Ind{{\mathbb{I}}}

\def\mes{{\rm mes\,}}

\def\const{{\rm const\,}}

\def\dist{{\rm dist\,}}
\def\R{{\bf R}}
\def\E{{\bf E }}
\def\P{{\bf P}}

\def\L{L}

\def\b{\beta}
\def\s{\d }
\def\g{\gamma}

\def\W{{\cal W}^*}

\def\oo{\bar}
\def\s{\sigma}

\def\p{\partial}
\def\G{\Gamma}
\def\GG{{\cal G}}

\def\L{{\cal L}}

\def\LLL{\w{\cal L}}

\newcommand{\be}{\begin{equation}}
\newcommand{\ee}{\end{equation}}
\newcommand{\bd}{\begin{displaymath}}
\newcommand{\ed}{\end{displaymath}}
\newcommand{\ba}{\begin{array}{ll}}
\newcommand{\ea}{\end{array}}
\newcommand{\baa}{\begin{eqnarray}}
\newcommand{\eaa}{\end{eqnarray}}
\newcommand{\baaa}{\begin{eqnarray*}}
\newcommand{\eaaa}{\end{eqnarray*}}   

\def\W{{\cal W}}

\def\ZZ{{\cal Z}}
\def\LLL{\w{\cal L}}

\def\CC{{\cal C}}

\def\G{D}

\def\AAA{{\mathfrak A}}
\def\W{\oo X}
\date{Submitted February 14, 2019. Revised: March 3, 2020 }
\title{Stochastic control problems and  HJB equations with excluded parameters of random inputs} 
\author{Nikolai Dokuchaev}
\begin{document}
\maketitle
\let\thefootnote\relax\footnote{School of Electrical Engineering, Computing and Mathematical Sciences, Curtin
University,  GPO Box U1987, Perth, 6845 Western Australia}
\begin{abstract}
This paper introduces a new type of second order stochastic  backward
Hamilton-Jacobi-Bellman (HJB) equations for optimal stochastic control problems
with a currently observable but  non-predicable parameter process,
in addition to the driving Brownian motion. The main feature of this
HJB equation is that it excludes   specifications of the parameter  process which dynamics
can be unspecified or unknown.
This allows to reduce the dimension of the state space. The paper considers  the case of control dependent diffusion coefficients and fully nonlinear HJB equations
under so-called Cordes conditions.

\end{abstract}
%

\section*{Introduction}
This paper considers optimal stochastic control problems in the continuous time setting.
The theory of these problems is well developed. In the diffusion Markovian setting, the value function
is usually represented by a parabolic  Hamilton-Jacobi-Bellman (HJB)  equation.
In non-Markovian  control problems,  the backward Hamilton-Jacobi-Bellman  equations
equations have to be replaced by corresponding backward SPDEs; this was first
observed  by Peng (1992).
\index{Pe1992}

Stochastic Partial Differential Equations (SPDEs)
 are well studied in the literature, including
 the case
 of  forward and backward  equations; see, e.g.,
\index{cite{Alos,Ba,DaPT,CK,CM,CMG, D92,D95,D05,D08b,D10,D12,DT, Duan,FZ, A1,Hu,kr,Ma,Ma99,Mas,Par,Roz,Wa,Zh}
 and the references
therein.} Walsh (1986),
Al\'os et al
(1999),
 Chojnowska-Michalik (1987),  Rozovsky (1990),  Zhou (1992),  Pardoux (1993), Bally {\it et al} (1994),  Chojnowska-Michalik and Goldys
 (1995),  Maslowski (1995), Da Prato and Tubaro (1996),
Gy\"ongy (1998),  Mattingly (1999),
  Duan {\em et al} (2003),  Caraballo {\em et al } (2004), Mohammed  {\em et al} (2008), Feng and Zhao (2012),
and the bibliography therein. In particular, backward SPDEs (BSPDEs) represent versions of
the so-called Bismut-Peng  equations where the diffusion
term is not given a priori but needs to be found; see e.g.   Hu and Peng (1991),  Peng
(1992),  Zhou (1992),  Dokuchaev (2008,2018)\index{(1992,1995,2008,2010, 2011,2012,2015a,b)}, Du  and Tang (2012), Du {\em at al} (2013), Hu and Peng (1995), Hu {\it et al} (2002),  Ma and
Yong (1999),  and the bibliography therein.   \index{The existing literature considers  BSPDEs of the
second order (and their generalizations) such that (i) the matrix of
the higher order coefficients is  positive, and (ii) the so-called
coercivity condition holds. }

In Bender and Dokuchaev (2016a,b) and Dokuchaev (2017), \index{\citet{BD1,BD2}} some special
BSPDEs
were derived for the value functions of special problems with linear state equations.
They represented  analogs of Hamilton-Jacobi-Bellman  equations for some non-Markovian
stochastic
optimal control problems associated with pricing of swing options in continuous time.
These equations are not exactly differential, since their solutions can be discontinuous in time, and they allow
 very mild conditions on the
underlying driving stochastic processes with unspecified dynamics.
 More precisely, the method does not have to assume a particular
evolution law of the underlying process;   the underlying processes
do not necessarily satisfy  stochastic differential equations of a known kind with a given structure.
In particular, the  First Order BSPDEs describe the value function
even in the situation where the underlying price process cannot be
described via a stochastic equation ever described in the
literature. The numerical solution requires just to calculate certain conditional expectations
of the functions of the process without using its evolution law (see the discussion in Section 4).
It can be also noted that these equations are not the same as the first order deterministic HJB equations known in the deterministic optimal
control.

The present paper extends these results on the setting with controlled diffusion with observed stochastic
parameter with unspecified dynamics.
The  paper considers a model where the controlled process is described as a stochastic
Ito process with coefficients depending on a random currently observable but unpredictable process $Z(t)$ being independent on the driving Brownian motion.
It is shown that the value function satisfies a second order BSPDE being a stochastic analog of the Hamilton-Jacobi-Bellman  equation
These stochastic equations are not exactly differential, since their solutions can be discontinuous in time, and they allow
 very mild conditions on the processes $Z(t)$  with unspecified dynamics. Similarly to the
 First Order BSPDEs introduced in  Bender and Dokuchaev (2016a,b), \index{\citet{BD1,BD2}}
the presented BSPDEs do not include  the parameters of a particular
evolution law of $Z(t)$;   these  processes $Z(t)$
do not necessarily satisfy  stochastic differential equations of a known kind with a given structure. This could be used to reduce the dimension of the equations even for the case when the dynamic of $Z(t)$ is known; see Remark \ref{remD} below.
The paper covers  the case of control dependent diffusion coefficients and fully nonlinear HJB equations
under so-called Cordes conditions.

The paper is organised as the following. In Section \ref{SecSet},
the control problem is described.  In Section \ref{SecB}, some background results
on weak solutions of Ito equations and related parabolic equations are provided.
 In Section \ref{SecM}, the main results on existence of optimal controls and backward
 SPDEs for the value functions are given. Section \ref{SecP}  contains the proofs.
\section{Stochastic control problem}\label{SecSet}
We are given an open domain $\G\subseteq\R^n$ such that either $\G=\R^n$
or $\G$ is bounded  with $C^{2+\a}$-smooth boundary $\p \G$ for some
$\a>0$; if $n=1$, then the condition of smoothness is not required.
Let $T>0$ be given, and let $Q\defi \G\times (0,T)$.
\par We are given a standard complete probability space $(\O,\F,\P)$ and  a $n$-dimensional Wiener process
$w(t)=(w_1(t),...,w_n(t))$, $t\ge 0$, with independent components such that $w(0)=0$.

Let an integer $d>0$ be given.
Let $Z(t)$, $t\ge 0$,  be a $d$-dimensional right continuous stochastic process with left limits
that is independent on $w$.

Let $\{\F_t^w\}_{t\ge 0}$ be the filtration generated by $w$, and
let $\{\F_t^{Z}\}_{t\ge 0}$ be the filtration generated by $Z$.

In addition, let $\{\GG_t\}_{t\ge 0}$ be the filtration generated by $(w,Z)$,
and
let $\{\oo\GG_t\}_{t\ge 0}$ be the filtration
such that $\oo\GG_t$ is the completion of all evens $\{A\cap B:\ A\in\F^w_t,\ B\in \F^Z_T\}$.

We assume that $\F^w=\F^Z_0=\GG_0$, and that this $\s$-algebra  is the completion of
the trivial $\s$-algebra  $\{\O,\emptyset\}$.

We denote by $\o$ the elements of the set
$\O=\{\o\}$.

Let  $\Delta \subset \R^m$  be a compact set.

Consider a  controlled Ito equations
\baa &&dy(t)=f(y(t),u(y(t),t),Z(t),t)dt+\beta  (y(t),u(y(t),t),Z(t),t)dw(t).\nonumber\label{3.1}
\eaa
For $t\ge s$ and a random vector $a$, we denote by $y^{a,s}(t)$,   the solution of this equation with the initial condition
\baa
&&y(s)=a.  \label{Ito}
\eaa
One selects functions    $u(\cdot):D \times[0,T]\times \O\to \Delta$ as controls.
\par Random vectors $a$ and  $y(t)$ take values in $\R^n$.

We assume that the  functions  $f(x,v,z,t): \R^n\times \Delta\times \R^d\times
\R \to \R^n$ and $\beta  (x,v,z,t):  \R^n \times\Delta\times\R^d\times  \R \to \R^{n\times n}$ are  continuous and bounded, together with the derivatives    $\frac{\p f}{\p  x}(x,v,z,t): \R^n\times \Delta\times \R^d\times
\R \to \R^{n\times n}$ and derivatives $\frac{\p^2 \b}{\p  x_k\p x_j}(x,v,z,t): \R^n\times \Delta\times \R^d\times
\R \to \R^n$ for $i,j=1,...,n$.

\par \Lindex{Consider  bounded functions  $\lambda(x,v,t): \R^n\times \Delta\times \R $    that are Borel measurable, continuous in   $v\in
\Delta$ for all  $x,t$.}

Consider a bounded Borel measurable
function  $\varphi(x,v,z,t): \R^n\times \Delta\times\R^d\times \R $ such that
the function $\varphi(x,\cdot,z,t):\Delta\to\R$ is continuous for all $x,z,t$.

\indexNN{NE NADO:Consider also  continuous functions  $\Phi(x): D\to \R$
such that $ \Phi|_{\p D}=0$. ????
Let $\Ind$ denotes the indicator function.}

\index{ For $s\in[0,T)$, let $\AAA_s$ be the set of all initial  vectors $a$,
such that $a\in D$ a.s., $a$ is $\GG_s$-measurable?,  and that
$a$ has a conditional probability distribution given $\GG_s$ that can described as the distribution}
 \par For $s\in[0,T)$, let $\AAA$ be the set of all initial random vectors $a$,
such that $a\in D$ a.s., $a$ is independent on $(w,Z)$,  and that
there exists $\rho\in H^{-1}$ such that
$(\rho,\Psi)_{H^0}=\E \Psi(a)$ for any $\Psi\in H^1$.
It can be noted that if $\rho\in L_1(D)$ then it is the probability density function of $a$. If $n=1$, then
the set $\AAA$ includes non-random $a\in D$; in this case, $\rho$ can be associated with the delta-function supported at $a\in D$.  In all cases, we say that $\rho\in H^{-1}$
describes the probability distribution of $a$.

For $a\in \AAA$ and some measurable function $\varphi:\R^n\times \Delta\times \R^d\times [0,T]\to\R$, let $\tau^{a,s}\defi \inf \{t:y^{a,s}(t)\notin D\}$, and let
\baa
&&F(a,u(\cdot))
\defi \index{\E \Ind_{\{\tau^{a,0}\ge T\}} \Phi
\bigl(y^{a,0}(T)\bigr)\Lindex{\exp\biggl\{-\int_0^{\tau^{a,0}}
 \lambda(y^{a,0}(r),r)dr\biggr\}}
+}\E\int_0^{\tau^{a,0}\land T}\varphi(y^{a,0}(t),u(t),Z(t),t)
\Lindex{\exp
\biggl\{-\int_s^t\lambda(y^{a,s}(r),r)dr\biggr\}}dt.
\eaa
\par
Consider control problem
\baa
\hbox{Minimize}\quad F(a,u(\cdot))\quad \hbox{over}\quad u(\cdot).
\label{P}
\eaa
\subsubsection*{Admissible controls}
Let us describe admissible controls.

Let  $U_0$
be the class of admissible control functions  $u(y,t,\o):\oo D \times [0,T]\times\O \to \Delta $ that
 are measurable and $\F^Z_t$-adapted for any $(y,t)\in \oo D\times [0,T]$.

\par
 Let $C(\Delta )$  be the space of real valued continuous functions
 defined on  $\Delta $. Let $C(\Delta )^*$  be its dual space.

 Let
$
\Delta_R=\bigl\{u\in C(\Delta )^*:\ \langle u, 1 \rangle_{C(\Delta)^*,C(\Delta)} =1,\
\langle  u,\nu    \rangle \ge 0\ ( \forall  \nu  \in C(\Delta ):\
\nu \ge 0)\bigr\}.
$
It can be noted that  $u\in \Delta_R$ can be associated with a probability measure on $\Delta$, so
$\langle  u,\nu    \rangle_{C(\Delta)^*,C(\Delta)}$  would be an averaging of the function $\nu(v)$ over $\Delta$ with resect to this
probability measure.

Let  $U_R$  be the class of admissible control functions  $u(y,t,\o):\oo D \times [0,T]\times\O\to \Delta_R$ that are  $\F^Z_t$-adapted for any $(y,t)\in \oo D\times [0,T]$.
\par
We assume that $U_0\subset U_R$, meaning that any $v\in \Delta$ is associated with the corresponding Dirac measure.

\section{Some background definitions and  results}
\label{SecB}
For a Banach space
$X$, we denote by $\|\cdot\|_{ X}$ the norm, and  we denote by $X^*$ its  dual space.
We will use notation $\langle a,b\rangle$ for $b(a)$, where $b:X\to \R$ is an element of $X^*$.

For a Hilbert space $X$, we denote by
 $(\cdot, \cdot )_{ X}$ the scalar product in  $X$.
 \par
 We denote Euclidean norm in $\R^k$ as $|\cdot|$, and $\bar G$ denotes
the closure of a region $G\subset\R^k$.
\par We introduce some spaces
of real valued functions.
\par
We denote by ${W_q^m}(D)$  the Sobolev  space of functions that
belong to $L_q(D)$ together with first $m$ derivatives, $q\ge 1$. In
particular, $$ \|u\|_{W_2^1(D)}\defi
\left(\|u\|_{L_2(D)}^2+\sum_{i=1}^n\left\|\frac{\p u}{\p
x_i}\right\|_{L_2(D)}^2 \right)^{1/2}.
$$
\par Let $H^0\defi L_2(D)$,
and let $H^1\defi \stackrel{\scriptscriptstyle 0}{W_2^1}(D)$ be the
closure in the ${W}_2^1(D)$-norm of the set of all smooth functions
$u:D\to\R$ such that  $u|_{\p D}\equiv 0$. Let $H^2=W^2_2(D)\cap
H^1$ be the space equipped with the norm of $W_2^2(D)$. The spaces
$H^k$ and $W_2^k(D)$ are called  Sobolev spaces, $k=0,1,2$; they are Hilbert
spaces, and $H^k$ is a closed subspace of $W_2^k(D)$, $k=1,2$.
As usual, we assume that $W_2^0(D)=H^0$.
\par
Let $H^{-k}\defi (H^k)^*$ be the dual spaces to the spaces $H^k$, $k=1,2$.

If $Y$ is the dual space for a space $X$,  we denote the dual pairing by
 $\langle\cdot ,\cdot \rangle_{Y,X}$.

\par
We denote by $\ell_k$ and $\oo\ell _{k}$ the Borel measure and the
Lebesgue measure in $\R^k$ respectively, and we denote by ${\cal
B}_{k}$ the $\sigma$-algebra of Borel sets in $\R^k$.  We denote
by $ \oo{{\cal B}}_{k}$ the completion of ${\cal B}_{k}$ with
respect to the measure $\ell_k$, or the $\sigma$-algebra of
Lebesgue sets in $\R^k$.
\par
We denote by $\oo{{\cal P}}^Z$  the completion (with respect to the
measure $\oo\ell_1\times\P$) of the $\s$-algebra of subsets of
$[0,T]\times\O$, generated by functions that are progressively
measurable with respect to $\F_t^Z$.
\par


For $k=-2,-1,0,1,2$, we  introduce spaces
 \baaa
&& X^{k}(s,t)\defi L^{2}\bigl([ s,t]\times\Omega,
{\oo{\cal P }}^Z,\oo\ell_{1}\times\P;  H^{k}\bigr),
\quad
\W^{k}(s,t) \defi L^{2}\bigl(\O, \F_T^Z,\P,L_2([s,t],\oo{{\cal B}}_1,\oo\ell_1,H^k\bigr),
\\ &&
\ZZ^k_t
\defi L^2\bigl(\Omega,\F^Z_t,\P;
H^k\bigr),\quad \CC^{k}(s,t)\defi C\bigl([s,t]; \ZZ^k_T\bigr).
\eaaa
The spaces $X^k(s,t)$, $\oo X^k(s,t)$, and $\ZZ_t^k$,  are Hilbert spaces.
\par
Further, we introduce spaces \baaa Y^{k}(s,t)\defi X^{k}(s,t)\!\cap
\CC^{k-1}(s,t), \quad \oo Y^{k}(s,t)\defi\oo X^{k}(s,t)\!\cap
\CC^{k-1}(s,t), \quad  k=1,2.
\eaaa
\index{with the norm $ \| u\| _{Y^k(s,t)}
\defi \| u\| _{{X}^k(s,t)} +\| u\| _{\CC^{k-1}(s,T)}$. $ \| u\| _{Y^k(s,t)}
\defi \| u\| _{{X}^k(s,t)} +\| u\| _{\CC^{k-1}(s,T)}. $}
\par
Note that spaces $X^k$ and $Y^k$ include adapted to $(w,Z)$ functions only, and the spaces $X^k$ and $Y^k$ include functions only that
are not necessarily adapted.

For brevity, we will use the notations
 $X^k\defi X^k(0,T)$,  $\oo X^k\defi\oo X^k(0,T)$, $\CC^k\defi \CC^k(0,T)$, $Y^k\defi Y^k(0,T)$,
and  $\oo Y^k\defi \oo Y^k(0,T)$.

\index{The same notations will be used for the spaces of vector and
matrix functions, meaning that all components belong to the
corresponding spaces.  In particular, $\|\cdot\|_{X^k}$ means the
sum of all this norms for all components.}

\subsubsection*{The assumption on the regularity of related linear parabolic equations}
The control problem (\ref{P}) is formulated for  a challenging case where the diffusion coefficients depend on the control.
This case is difficult even in the Markovian  setting (i.e. where $Z(t)\equiv 0$), because
 corresponding forward and backward Kolmogorov parabolic equations for distributions
are  equations in non-divergent form with
discontinuous coefficients at higher derivatives; they
do not feature sufficient regularity in a case of non-smooth closed loop controls $u(\cdot,t)$.  Their
investigation is most complicated because, in general, in the case
of discontinuous coefficients, the uniqueness of a solution for
nonlinear parabolic or elliptic equations can fail, and there is
no a priori estimate for partial derivatives of a solution; see. e.g. Krylov  (1987) and the literature therein.
On the other hand, a typical optimal control is not expected to be smooth.

There are two main approaches to overcoming these difficulties via relaxation of the requirements for the solution.
One approach is to consider the so-called viscosity solutions; see, e.g.,  Crandall and Lions
(1983). \index{CL}  Another  approach is to accept  solutions with measure-valued  second derivatives; see e.g.  Krylov (1980,1987). \index{KR,KR2}
In this paper, we will not be using either of these approaches. Instead, we restrict our consideration
only by the cases where the backward Kolmogorov equations
for controls $u\in U_R$ features solutions with  $L_2$-integrable derivatives. This still can cover some important case. For example, this setting covers the case where  the part of the diffusion coefficient depending on the control is restricted in size (Condition \ref{condA} below).

 Let $b(x,v,z,t)\defi \beta(x,v,z,t) \beta(x,v,z,t)^\top /2$. We assume that
 \baaa
 \inf_{\xi\in \R^n, (x,t)\in Q, v\in \Delta, z\in\R^d}\frac{\xi^\top b(x,v,z,t) \xi}{|\xi|^2}>0.
 \label{delta}\eaaa

We assume that the domains for  $b, f,\Lindex{ \lambda} $  and $\varphi $ are extended
to   $\R^n\times \Delta_R\times \R$ as the following. For
$(x,u,z,t)\in   \R^n\times \Delta_R\times \R^d\times \R $, we assume that
$f(x,u,z,t)=\langle f(x,.,z,t),u \rangle_{C(\Delta)^*,C(\Delta)}$ and $b(x,u,z,t)=\langle
b(x,.,z,t),u \rangle_{C(\Delta)^*,C(\Delta)} $.    \Lindex{etc.}  It can be noted that  $u\in \Delta_R$ can be associated with a probability measure on $\Delta$; therefore,  these extensions represent averaging over this measure.
\par

 Let us define differential operators
\baa
A(x,u,z,t){\rm V}=\sum_{i,j=1}^n  b_{ij}(x,u,z,t)\frac{\p^2{\rm V}}{\p  x_i \p
x_j} (x)+ \sum_{i=1}^n f_i(x,u,z,t)\frac{\p  {\rm V}}{\p x_i}(x)\comm{ - \lambda
(x,u,t) {\rm V}(x)},  \label{3.4}
 \eaa where $b_{ij}$, $f_i,  x_j$  are components
of $f,x$, $b$. Since these functions  depend on $Z(t)$, i.e., these
operators have  random coefficients. \par
  In $Q=D\times
(0,T)$, consider, for  an admissible  $u(\cdot)$, a boundary value problem
\baa
&&\frac{\p {\rm V}}{\p  t}(x,t)+A(x,u(x,t),Z(t),t){\rm V}(x,t)=-\psi
(x,t),  \label{parab0}\nonumber
\\ && {\rm V}|_{x\in \p D}=0, \quad
{\rm V}(x,T)=0 \index{\Psi(x)}.  \label{parab}
 \eaa
\begin{condition}\label{condA}  Assume that the function  $b$ is such that,
for any $u\in U_R$,
the problem   (\ref{parab})  has a unique solution   ${\rm V}\in
\oo Y^2$  for any $\psi\in X^0$.\index{and  $\Psi\in H^1$.} Moreover,
\baa  \int_0^T\left( \left\| \frac{\p {\rm V}}{\p t}(\cdot,t) \right\|_{H^0}^2+\|{\rm V}(\cdot,t) \|_{H^2}^2\right)dt \le   c \|\psi\|_{L_2(Q)}^2 \index{ + \|\Psi \|_{H^1}^2)}\quad \hbox{a.s.}, \eaa where
$c>0$ depends on $D,n,T,f,b\Lindex{,\lambda}$ only.
\end{condition}

The following result  describes some
special cases where Condition \ref{condA} holds.
\begin{lemma} \label{lemmaA}
Condition \ref{condA} holds if at least one of the following conditions is satisfied.
\begin{enumerate}
\item  $b(x,v,z,t)\equiv b(x,z,t)$;
\item The matrix  $b$ has the
form $b(x,v,z,t)=\oo b(x,z,t)+\w b(x,v,z,t), $ where $\oo b(x,z,t)=\oo
b(x,z,t)^\top$ is a continuous bounded matrix such that
 \baa
C_b\defi  \inf_{\xi\in \R^n, (x,t)\in Q, z\in\R^d}\frac{\xi^\top\oo b(x,z,t) \xi}{|\xi|^2}>0,
 \label{delta2}\eaa
and where $$
\sup_{(x,t)\in  Q,\ v\in\Delta, \ z\in\R^d}\sum^n_{i,k=1}\w b_{ik}(x,v,z,t)^2
<\frac{C_b^2}{n}.
$$
\item  The matrix  $b$ has the
form $b(x,v,z,t)=\oo b(x,z,t)+\w b(x,v,z,t), $ where $\oo b(x,z,t)=\oo
b(x,z,t)^\top$ is a continuous bounded matrix such that (\ref{delta}) holds.
The matrix function $\w b(x,v,t)$ is
symmetric and such that there exists a set ${\cal
N}\subseteq\{1,\ldots,n\}$ such that
$$\w b_{ij}\equiv \w b_{ji}\equiv 0 \quad \forall i,j: i\notin
{\cal N}, j\notin  {\cal  N},
$$ and there exists a set $\{\gamma _k\}_{k\in {\cal  N}}$ such that $\g_k\in(0,2)$ for all $k$
and
\baaa
\left(\sum_{k\in\N}\frac{1}{2\gamma_k}\right)
\sup_{(x,t)\in  Q,\ v\in\Delta, \ z\in\R^d}\sum_{k\in \N} \Biggl(\sum_{i\in\N}
\w b_{ik}(x,v,z,t)^2  &+&\!4\sum_{i\notin\N}\w b_{ik}(x,v,z,t)^2
\\
&+&\frac{\gamma_k}{2-\gamma_k}\w b_{kk}(x,v,z,t)^2\Biggr) <\! C_b^2.
\eaaa
\end{enumerate}
\end{lemma}

The result of Lemma \ref{lemmaA} was obtained  in  Dokuchaev (1997); a related result can be found in Dokuchaev (2005).

Conditions (ii) and (iii)  represents  analogs of the
so-called Cordes conditions that ensure regularity of solutions of
boundary value problems for second order equations and that are
known as Cordes conditions.

\subsubsection{On Cordes conditions: some historical remarks}
It is known that discontinuity of the higher order coefficients  for linear parabolic and elliptic equations
in non-divergent form causes problems with regularity of the solutions.
This makes analysis of corresponding HJB equations difficult; see, e.g. Krylov  (1987) and the literature therein.
A possible approach is to consider so-called viscosity solutions; see, e.g., Fleming and Soner (1993).
These solutions may not have all integrable derivatives.
In some cases, it is still possible to have solutions with $L_2$-integrable derivatives;
  this can be achieved with  some restrictions the scale of discontinuities
for the higher order  coefficients.
The original  Cordes  conditions
restricts the scattering  of the eigenvalues  of the matrix of the
coefficients at higher derivatives (see Cordes (1956)).  Related conditions from
 Talenti (1965), Koshelev (1982), Kalita (1989),  Landis (1998), on
the eigenvalues are also called Cordes type conditions.   A closed
condition is presented  implicitly in  the proof of the uniqueness of a weak
solution  in  Gihman
and Skorohod (1975), Section 3 of Chapter 3.

Cordes (1956) considered elliptic equations.   Landis
(1998) considered both elliptic and parabolic equations.  Koshelev
(1982)  considered systems of elliptic equations of divergent type
and H\"older property of solutions.  Kalita (1989) considered
union of divergent and nondevirgent cases.
\par
Conditions from  Cordes (1956) are such that they are not
necessary satisfied even for constant non-degenerate matrices $b$,
therefore, the condition for $b=b(x)$ means   that the
corresponding inequalities are satisfied for all $x_0$ for some
non-degenerate matrix $\theta(x_0)$  and
 $\ww  b(x)=\theta (x_0)^T  b(x) \theta (x_0)$,
where    $x$  is from $\e$-neighborhood of  $x_0$  ($\e>0$ is
given). Conditions  in Lemma \ref{lemmaA}  ensure  also solvability and uniqueness for first boundary
value problem for nondivirgent parabolic equation with
discontinuous diffusion coefficients. Moreover, conditions in Lemma  \ref{lemmaA}
 ensure prior estimate required in Condition \ref{condA}, in contrast
with the existing  literature.  Second order SPDEs satisfying Cordes conditions similar to the ones
in Lemma \ref{lemmaA} were considered in  Dokuchaev (2005) for forward stochastic SPDEs.
Some comparison of different types of
Cordes conditions can be found in Dokuchaev (1997,2005).

\subsubsection*{Some auxiliary operators}
 It can be seen that Condition \ref{condA} implies continuity,   for any $u\in U_R$, of the
 following linear operators
 \baaa
&&\oo L(u(\cdot)): X^0\to   \W^2, \indexNN{\qquad \oo {\mathfrak L}(u(\cdot)):X^0\to   \ZZ_T^1,}
 \quad L(u(\cdot)): X^0\to   X^2, \qquad {\mathfrak L}(u(\cdot)):X^0\to  H^1
 \eaaa
 defined such that
 \baaa
&& \oo V=\oo L(u(\cdot))\psi,\qquad \indexNN{\oo V(\cdot  ,0)=\oo{\mathfrak L}(u(\cdot))\psi,\quad} V= L(u(\cdot))\psi,\qquad V(\cdot  ,0)={\mathfrak L}(u(\cdot))\psi,\eaaa
where $\oo V$ is the solution of the problem (\ref{parab}), \index{with $\Psi=0$,} and where
 $V(\cdot,t)=\E\{\oo V(\cdot,t)|\F_t^Z\}$. (Remind that $\ZZ_0^k=H^k$).
The corresponding  adjoint operators \indexNN{NE NADO $L(u(\cdot))^*:X^{-2}\to X^0$}
 ${\mathfrak L}(u(\cdot))^*:H^{-1}\to X^0$ are linear and continuous as well.

\indexNN{NE NADO  In addition,  we introduce operators $L_{s,T}:  L_2(D\times
(s,t))  \to H^1$,
${\cal     L}_{s,T}:H^1     \to      H^1$  such that
$V(\cdot ,s)=L_{s,T}\psi+\LLL_{s,T}\Psi $ for
$V=L (u(\cdot))\psi+{\mathfrak L}(u(\cdot))\Psi $.
 By Condition \ref{condA}, these linear
operators are continuous. The corresponding adjoint operators
$L^*_{s,t}:H^{-1} \to L_2(D\times [s,t])$ and
$\LLL^*_{s,t}:H^{-1} \to H^{-1}$ are linear and continuous as well.
}
 \subsubsection*{On weak solutions of Ito equations}
We consider solutions of  (\ref{Ito}) for $u\in U_R$. In this case, we assume that
$\b$ is  defined  as a square-root of the matrix
$2b$; to ensure that the choice of the square-root is unique,
we can require, for example, that  the matrix $\b$ is positive-definite everywhere.
\par
It can be noted that $b(\cdot,u,\cdot)$ is affine in $u\in U_R$, but this
is not necessarily  the case for   $\b(\cdot,u,\cdot)$.

We consider weak solutions of  (\ref{Ito}) for $u\in U_R$  as described  in the following lemma.
\begin{lemma}\label{ThEI}   Let  $a \in L^2(\Omega,\F,\P,\R^n)$ be independent on $w$. Then, for any $u\in U_R$,  there exists
a set
$$
\biggl\{  (\w\O,\w\F,\w\P),  (w(t),{\cal   F}_t),
y^{a,0}(t) \biggr\},
$$
where    $(\w\O,\w\F,\w\P)$ is a probability space  such that
 $a \in L^2 (\w\Omega,\F,\P)$,
$(w(t),\F_t)$  is a $n$-dimensional Wiener process on
$(\w\O,\w\F,\w\P)$,  $\F_t\subseteq  \w\F$ is a
filtration of $\s$-algebras of events such that  $w(t)-w(s)$
do not depend on $(a,Z)$ and on $\F_s$ for  $t>s$, and
$y^{a,0}(t)$ is the solution of (\ref{Ito}) for this $w(t)$.
\end{lemma}
\par
This theorem can be found, in particular, in Krylov (1980), Chapter 2;
it is formulated therein
for non-random $(a,Z(\cdot))$, which is
unessential  since $(a,Z)$ are independent on $w$.

The results of  Lemma \ref{ThWE}-\ref{ThPDF} below were obtained in Dokuchaev (1997).

\begin{lemma}\label{ThWE} Let  $a\in\AAA$,  and let
$\rho$ describes the probability distribution of $a$.  Assume that Condition \ref{condA} is satisfied. Then, for any $u(\cdot)\in U_R$, equation (\ref{Ito})  has a unique
weak solution, meaning that the solution  is univalent
with respect to the probability distribution.
\end{lemma}

In this paper, we consider problem (\ref{P}) with weak solutions of
equation (\ref{Ito}).
 \begin{lemma}\label{ThEIF}  
 Let  $a\in\AAA$,  and let
$\rho$ describes the probability distribution of $a$.
Then, for any measurable $u\in U_R$,
 \baaa
\E\left\{\int_0^{\tau^{a,0}\land T}\varphi(y^{a,0}(t),u(t),Z(r),t)
\Lindex{\exp
\biggl\{-\int_0^t\lambda(y^{a,s}(r),Z(r),r)dr\biggr\}}dt\Biggl|\F_T^Z\right\}
=\bigl\langle\rho, \oo V(\cdot ,0)\bigr\rangle_{H^{-1},H^1}\quad\hbox{a.s.}
\eaaa
and
\baaa
 F(a,u(\cdot))=\E \bigl\langle\rho, \oo V(\cdot ,0)\bigr\rangle_{H^{-1},H^1}
 =\bigl\langle\rho,  V(\cdot ,0)\bigr\rangle_{H^{-1},H^1}
 \label{3.4.2'}
\eaaa
where  $\oo V=\oo L(u(\cdot))\varphi(\cdot,u(\cdot),Z(\cdot),\cdot)\index{+\L(u(\cdot))\Phi}\in\oo Y^2$, and $V=L(u(\cdot))\varphi(\cdot,u(\cdot),Z(\cdot),\cdot)\index{+\L(u(\cdot))\Phi}\in Y^2$. In addition,
$$
 |F(a,u(\cdot))|\le  c \|\rho\|_{H^{-1}}\|\varphi(\cdot,u(\cdot),Z(\cdot),\cdot)\|_{X^0}\index{+ \|\Phi \|_{H^1}},
$$
where  $C>0$ is a constant occurring in Condition
\ref{condA}.
\end{lemma}
By Lemma \ref{ThEIF}, the functionals  $F(a,u(\cdot))$ are defined for Borel measurable $u(\cdot)\in  U_R$, and   $F(a,u(\cdot))=\bigl\langle\rho,  V(\cdot ,0)\bigr\rangle_{H^{-1},H^1}$, where $V=L(u(\cdot))\varphi (\cdot ,u(\cdot),Z(\cdot),\cdot)$. However,   the value $\bigl\langle\rho,  V(\cdot ,0),\rho\bigr\rangle_{H^{-1},H^1}$  is defined also   for $u(\cdot)\in
U_R$, and it does not depend on the choice of a representative  of a class
of equivalency;  for  $u(\cdot)\in   U_R$
there exists a Borel measurable equivalent function $u(\cdot)\in  U_R$.
Respectively, we presume  that the functionals $F(a,u(\cdot))$
are extended on  $u(\cdot)\in U_R$.

 \begin{lemma}\label{ThPDF} 
 Under the assumptions of Lemma \ref{ThEIF}, the weak solution $y^{a,0}(t)$  of equation  (\ref{Ito}) with $s=0$, considered \Lindex{  being killed inside $D$ with the rate $\lambda$ and being killed} on the boundary $D$,  has the conditional distribution given
$\F^Z_T$ featuring the probability density
function $p=\L(u(\cdot))^*\rho\in L_2(Q)$ a.s.. Moreover, $p \in X^0$.
\end{lemma}
%

\begin{corollary}\label{corrMP} (The Maximum Principle). Assume
that conditions of Lemma  \ref{ThPDF}  are satisfied and, in
addition, that
 $\varphi(x,v,z,t)\ge 0$ for a.e. $x,t$ for all $v\in \Delta_R$, $z\in\R^d$,
\indexNN{ $\Phi (x)\ge 0$,}  $\langle\rho,\phi\rangle_{H^{-1},H^1}\ge 0$ a.s. for all $\phi\in H^1$ such that $\phi(x)\ge 0$ a.s.., where $\langle\cdot,\cdot\rangle_{H^{-1},H^1}$ denote the natural pairing between $H^{-1}$ and  $H^1$. Then, for any $u\in U_R$, we have that  $V(x,t,\o)\ge 0$ and
 $p(x,t,\o)\ge 0$ a.e., \indexNN{all $t$ for a.e. $x,t$ a.s.,} where $V=L(u(\cdot))\varphi+\L(u(\cdot)) \Phi$ and
    $p=\L(u(\cdot))^*\rho$.
    \end{corollary}
 The following lemma provides a strengthened   version of the maximum principle.
 \begin{lemma}\label{ThNN} Let  $a\in\AAA$,  and let
$\rho\in H^{-1}$ describing  the probability distribution of $a$ be  such that
 $\langle \rho,\phi\rangle_{H^{-1},H^1}>0$ for all $\psi\in H^1$ such that  $\mes\{x\in D:\psi(x)>0\}$>0.
 Then, under the assumptions of Lemma
\ref{ThPDF}, $p(x,t,\o)>0$ for a.e.. $x\in D$, $t\in(0,T]$, $\o\in\O$
\end{lemma}

\section{The main results}\label{SecM}
Up to the end of this paper, we assume  that Condition \ref{condA} holds.
\indexNN{NE NADO:
Let  $U^s_R$
be the set   $U_R$   provided with the topology
of the dual space $L^1(\Theta
\times [0,T]\to C(\Delta  ))$, i.e., topology of   $L^{\infty}(Q\to
C^*_s(\Delta  ))$.
 \begin{theorem}\label{Th3.1} Let  $a\in \AAA$. Then
$F(a,u(\cdot))$ is continuous on $U_R^s$.
\end{theorem}
}
\begin{theorem}\label{ThEr} Let  $a\in \AAA$. Then
there  exists an optimal solution  $u\in U_R$ for problem (\ref{P}).
\end{theorem}


\begin{theorem}\label{ThM} Assume that the set $\Delta$ is convex. In this case, the following holds.
\begin{enumerate}
\item
 There exists  $\w u(\cdot)\in U_0$
such that $\w V=L(\w u(\cdot))\varphi(\cdot,\w u(\cdot),Z(\cdot),\cdot)$
 satisfies  the following modification of the Hamilton-Jacobi-Bellman equation
  \baa
  &&\w V(x,t)=\E\left\{ \int^T_t\inf_{v\in \Delta}\left[
  A(x,v,Z(s),s)\w V(x,s)+\varphi (x,v,Z(s),s)\right]ds\Biggl| \F_t^Z\right\},  \label{3.7}\nonumber \\ &&\w V(x,t) |_{x\in \p D}=0, \quad \w
V(x,T)=0,  \label{HJB}
\eaa
\index{ \baa
  &&\frac{\p
\w V}{\p      t}(x,t)+A(x,\w u(x,t),Z(t),t)\w V(x,t)+\varphi (x,\w
u(x,t),t)=0,  \label{3.7}\\ &&\w V |_{x\in \p D}=0, \quad \w
V(x,T)=0,  \label{HJB}  \\&&A(x,\w u(x,t),Z(t),t)\w V(x,t)+\varphi (x, \w
u(x,t),Z(t),t)\nonumber\\ &&\hphantom{xxxxxxxxxxxx}\le  A(x,v,Z(t),t)\w V(x,t)+\varphi (x,v,Z(t),t) \qquad  \label{3.9}
\eaa}
for all $t\in[0,T|$ for a.e.    $x\in  D$ a.s..

\item  If  $\w u\in U_0$ is such as described above, then $F(a,\w u(\cdot))=\bigl\langle\rho,  \w V(\cdot ,0)\bigr\rangle_{H^{-1},H^1}$ and
\baa F(a,\w u(\cdot))\le F (a,u(\cdot)) \quad
 \forall  u(\cdot)\in U_0, a\in \AAA.
\label{3.11}\eaa
\end{enumerate}
\end{theorem}

It can be noted that  Lemma \ref{lemmaA} ensures that $\w V\in Y^2$; this  implies a regularity of the solution of (\ref{HJB}).
 \par
\begin{remark}\label{remD} {\rm Equation (\ref{HJB}) does not include  the parameters of a particular
evolution law of $Z(t)$;   these  processes $Z(t)$
do not necessarily satisfy  stochastic differential equations or jump-diffusion equation of a known kind with a given structure.
This could be used to reduce the dimension of the equations even for the case when the dynamic of $Z(t)$ is known.
Assume, for example,
that the controlled process $y(t)$ is $n$-dimensional and
that  a scalar process $Z(t)$   is defined as $Z(t)=C X(t)$, where  $C\in\R^{1\times N}$,
and where $X(t)$ is a $N$-dimensional  solution of an Ito equation. Then the traditional parabolic
HJB equation in Markovian setting   would require the state space $\R^{n+N}$. On the other hand,  the state space for equation (\ref{HJB}) is $\R^n$; this gives a significant  dimension reduction for large $N$. }
\end{remark}
\section{Proofs}\label{SecP}
\subsection{Proof of Lemma \ref{ThNN}}
Let $\P_0\{\cdot|(a,Z)\}$ be a probability measure that is equivalent to $\P\{\cdot|(a,Z)\}$
and such that the process $y^{a,0}(t)$ is a martingale on a the conditional probability
space given $(a,Z)$; this measure exists by Girsanov Theorem. In  this case,
for any $\a\in D$ and any $\e>0$,
\baaa
\P_0(\sup_{t\in[0,T]} |y^{\a,0}(t)-\a|\le \e\,|\,\F_T^Z\}>0\quad \hbox{a.s.}.
\label{aD1}\eaaa
This follows from the properties of standard one-dimensional
Wiener processes and from  the Dambis-Dubins-Schwartz theorem applied to the components of
the vector process  $y^{a,0}(\cdot)$.

 In  this case,
for any $\a\in D$ and any $\e>0$,
\baaa
\P_0(\sup_{t\in[0,T]} |y^{\a,0}(t)-\a|\le \e\,|\,\F_T^Z\}>0\quad \hbox{a.s.}.
\label{aD2}\eaaa
This implies that
\baa
\P(\sup_{t\in[0,T]} |y^{\a,0}(t)-\a|\le \e\,|\,\F_T^Z\}>0\quad \hbox{a.s.}.
\label{aD3}\eaa

Further, suppose that there exists  a
domain $D_0\subset D$ such that
$\P(y^{a,0}(T)\in D_0\,|\,\F^Z_T)=0$. Let $D_\e\defi \{x\in D_0:\  \dist(x,\p D_0)>\e\}$.
Let $\e>0$ be such that $\mes D_\e>0$.

By the properties of $\rho$, it follows  that $\P(a\in D_\e)>0$.
 Clearly, (\ref{aD3}) implies that  \baaa
\P(y^{a,0}(t)\in D_0\quad \forall t\in[0,T]\,|\,\F^Z_T,\ a\in D_\e)>0\quad \hbox{a.s.}.
\label{aD4}\eaaa
This completes the of Lemma \ref{ThNN}.
$\Box$

\subsection{Proof of Theorem \ref{ThEr}}
We denote by  $U_R^w$ the set  $U_R$ provided
with the weak topology of the space being dual to the space
$L^1   (D  \times [0,T]\to  C(\Delta ))$. It is known
 that  the set  $U_R^w$  is convex, compact, and sequentially compact; see e. g. Varga (1972), Ch.IV.

 For $u\in U_R$, the operators
 $ A(\cdot,u(\cdot),\cdot,t)\defi A(x,u(x,t),Z(t),t): \ZZ_T^2\to \ZZ_T^{0}$ are continuous, hence the adjoint  operators
$A^*(\cdot,u(\cdot),\cdot,t): \ZZ_T^0\to \ZZ_T^{-2}$ are continuous; here  $\ZZ_T^{-2}= (\ZZ_T^{2})^*$.
These operators are such that
\baaa
\langle  A^*(\cdot,u(\cdot),\cdot,t)\eta,\psi\rangle_{ \ZZ_T^{-2}, \ZZ_T^{0} } =(\eta, A(\cdot,u(\cdot),\cdot,t) \psi)_{\ZZ_T^0},\quad \psi\in \ZZ_T^2,\ \eta\in \ZZ_T^0,\quad t\in[0,T].
\eaaa

%

 \begin{proposition}\label{prop3.2} Let $\rho\in H^{0}$ describes the probability distribution of $a$
 (in particular, if $\rho\in H^0$ then $\rho$ is the density for $a$).
 Let $u_i(\cdot)\in U_R$ be such that the derivatives
up to the second order for $u_i(x,t,\o)$ with respect to $x$ are bounded,
$\alpha_i\ge 0$, $i=1,...,N$, $\sum_{i=1}^N \alpha_i=1$, $N=1,2,3,...$.
Let  $ p_i =\L(u_i(\cdot))^*\rho$.  Let $\ww p(x,t)=  \sum_{i=1}^N  \alpha_i
p_i(x,t)$. Let us consider control  $$\ww u(x,t)=\ww p(x,t)^{-1}
\sum_{i=1}^N \alpha_i p_i(x,t)u_i(x,t). $$ (It can be seen that
$\ww  u\in  U_R$ since $\Delta_R$ is a convex set). Let  $p_{\ww u}=\L^*(\ww u(\cdot))\rho
$. Then \baaa
p_{\ww u}=\ww p,\qquad \sum_{i=1}^N \alpha_i
F(a,u_i(\cdot))=F(a,\ww u(\cdot)).
\eaaa
\end{proposition}
\par

{\em Proof of  Proposition \ref{prop3.2}}.
In this proof, we use that the coefficients $f$ and $b$ defined as $f(x,v,z,t)=\langle f(x,\cdot,z,t),v\rangle_{C^{*}(\Delta),C(\Delta)} $ and
$b(x,v,z,t)=\langle b(x,\cdot,z,t),v\rangle _{C^{*}(\Delta),C(\Delta)} $ are affine with respect to the probability measures
$v\in \Delta_R$.  It can be noted that, however, that it is not required that the functions
$f(x,\cdot,z,t):\Delta\to\R^n$ and $b(x,\cdot,z,t):\Delta\to\R^{n\times n}$ are affine with respect to $u\in \Delta$.

Since the function    $\varphi(x,v,Z(t),t)$
\Lindex{$\lambda (x,u,t)$,  $f(x,v,Z(t),t)$,} is affine in  $v\in  U_R$, we have that
\baaa
\sum_{i=1}^N \alpha_i F(a,u_i(\cdot))=
\E\int_Q dx\,dt\sum_{i=1}^N \alpha_i p_i(x,t)\varphi
(x,u_i(x,t),Z(t),t)\\= \E\int_Q \w p(x,t) \varphi_i(x,\ww
u(x,t),Z(t),t)dx\,dt  . \eaaa
Therefore, it suffices to show that
$p_{\ww u}=\ww p$.

Let us assume first that $\rho\in H^0$.
Since we selected smooth in $x$ controls $u_i=u_i(x,t,\o)$, we have that
 satisfy the forward Kolmogorov equations
\baa
  &&\frac{\p p_i}{\p      t}(x,t)=A^*(x,u_i(x,t),Z(t),t)p_i(x,t),  \nonumber
  \\ &&p_i(x,t,\o)|_{x\in \p D}=0, \quad \w
p_i(x,0,\o)=\rho(x). \label{pi}
\eaa
Here $A^*(\cdot ,u_i(\cdot,t),Z(t),t):H^1\to H^{-1}$ are the formally operators
for the operators $A(\cdot ,u_i(\cdot,t),Z(t),t):H^1\to H^{-1}$.

The classical theory for these parabolic equations
ensures that these equation have unique solutions  $p_i\in Y^1$; moreover, we have that $p_i\in Y^2$ as well; see, e.g., \index{la}  Ladyzhenskaya (1985),  Sections III.4-III.5).

  We sum  Kolmogorov's equations for
$p_i$  (or, more precisely, for     $ \alpha_i p_i$), using the fact that the functions   $b(x,v,Z(t),t)$,    and $f(x,v,Z(t),t)$   \Lindex{$\lambda (x,v,Z(t),t)$,}
   are affine in $u\in U_R$, meaning relations such as
\baaa
\sum_{i=1}^N \alpha_i p_i(x,t) f(x,u_i(x,t),Z(t),t) =\langle
f(x,\cdot,t),\sum_{i=1}^N \alpha_i p_i(x,t)u_i(x,t) \rangle_{C^{*}(\Delta),C(\Delta)}\\
=\ww p(x,t) \langle f(x,\cdot,Z(t),t),u(x,t) \rangle _{C^{*}(\Delta),C(\Delta)}=\ww p(x,t)f(x,\ww u(x,t),t). \eaaa

It gives that
\baaa
\sum_{i=1}^N \alpha_i  A^*(\cdot  ,u_i(\cdot),\cdot)p_{u_i}(\cdot ,\cdot) =A^*(\cdot ,\ww u(\cdot),\cdot)p_{\ww
u}(\cdot ,\cdot);
\eaaa
the equality here is in  $X^{-2}$. From the form of the corresponding forward
Kolmogorov equation for the conditional density $p$ on the conditional probability space
given $Z$, we obtain that  $p_{\ww u}=\ww p$ if $\rho\in H^0$.
Since $H^0$ is everywhere dense in $H^{-1}$, it follows that this identity holds also in the case
if $\rho\in H^{-1}$.
This completes the proof of  Proposition \ref{prop3.2}. $\Box$

\par
Further, let $J\defi \inf_{U_R}  F(a,u(\cdot))$.  Clearly, there exists
a sequence of controls  $u_i(\cdot)\in U_R$,
$i=1,2,...$,    such that the derivatives
up to the second order for $u(x,t,\o)$ with respect to $x$ are bounded,   and that
 $F(a,u_i(\cdot))\to  J$  as $i\to +\infty $.
Let  $\oo p_i\defi \L(u_i(\cdot))^*\rho $.  By passing to a
weakly converging
subsequence, we see that  there exists  $\w p\in X^0$ such that
 $\oo p_i \to
\w p$  as  $i\to +\infty $ weakly in $X^0$.  By  the  Mazur Theorem  (see, e.g., Yosida (1995),  p.173), there
exists a sequence of convex combinations $p_i(\cdot)=\sum_{j=1}^i    \alpha_j
\oo p_j(\cdot)$, $\alpha_j=\alpha_j(i)$, $\alpha_j\ge
0$, $\sum_{j=1}^i\alpha_j=1$ such that  $p_i(\cdot) \to  \w p(\cdot
)$   in $X^0$. By Proposition \ref{prop3.2},  for any
$p_i(x,t)$ there exists  $v_i(\cdot)\in U_R$, such that
$p_i(x,t)=\L^*(v_i(\cdot))\rho$ and  $F(a,v_i(\cdot))=\sum_{j=1}^i \alpha_j F(a,u_j(\cdot))$.
\indexNN{NE NADO In this case, $F(a,u(\cdot))=(p,\varphi  (\cdot,u(\cdot),\cdot))_{X^0}+(p(\cdot,T),\Phi)_{\ZZ_T^0}$, where  $p=\L(u(\cdot))^*\rho $.}  In addition, it is easy to see that
$F(a,v_i(\cdot))\to J$ as  $i\to +\infty$.

The set $U_R$ is compact in the topology of   $U_R^w$.  Passing to a
subsequence, we obtain that there exists  $\w u(\cdot)$ such that $v_i(\cdot) \to  \w u(\cdot)$
as   $i\to +\infty $ in the topology of   $U_R^w$.

Let us define operators  \baaa
A_i(\cdot,t)\defi A(x,v_i(x,t),Z(t),t): \ZZ_T^2\to \ZZ_T^{0},\quad \w A^*(\cdot,t)\defi A^*(x,\w u(x,t),Z(t),t): \ZZ_T^0\to \ZZ_T^{-2}.\eaaa These operators are continuous.

By the assumptions on $u_i$ and $v_i$, it follows that $p_i$ belong to $Y^1$ and that $p_(t)$ represent conditional densities
given $\F_t^Z$ for processes $y^{a,0}(t)$  being killed on the boundary. Hence they
satisfy the forward Kolmogorov equations
\baaa
p_i(\cdot,t)=\rho+\int_0^t A_i^*(\cdot,s)p_i(\cdot,s)ds.
\eaaa

Furthermore,
\baaa
A_i^*p_i- \w A^*\w p=R_{1,i}+R_{2,i},
\eaaa
where
\baaa
&&R_{1,i}\defi A^*_i\w p-\w A^*\w p,\qquad R_{2,i}\defi  A^*_ip_i  -A^*_i\w p,
\eaaa
i.e. $R_{1,i}= [A^*_i-\w A^*]\w p$ and  $R_{2,i}= A^*_i[p_i  -\w p]$.

Clearly, $\|R_{2,i}\|_{X^{-2}}\to 0$ as $i\to +\infty$ and
\baaa
R_{1,i} \to 0 \quad\hbox{weakly in}\quad X^{-2} \ \ \ \hbox{as}\ \ i \to +\infty .
\eaaa
Hence, for any $t\in[0,T]$,
\baaa
\int_0^t A^*_i(\cdot,s)p_i(\cdot,s)ds\to  \int_0^t \w A^*(\cdot,s)\w p(\cdot,s)ds\quad\hbox{weakly in}\quad \ZZ_T^{-2} \ \ \ \hbox{as}\ \ i \to +\infty.
\eaaa
It follows that
\baaa
\w p(\cdot,t)=\rho+\int_0^t \w A^*(\cdot,s)\w p(\cdot,s)ds.
\eaaa
These equalities hold for all $t\in[0,T]$ in $\ZZ_T^{-2}$.

Similarly, we obtain that, for any $\psi\in \ZZ_T^2$,
\baaa
\langle  \w p(\cdot,t),\psi\rangle_{ \ZZ_T^{-2}, \ZZ_T^{0} }
&=&\langle  \rho,\psi\rangle_{ \ZZ_T^{-2}, \ZZ_T^{0} } +\int_0^t
\langle \w  A^*(\cdot,s)\w p(\cdot,s),\psi\rangle_{ \ZZ_T^{-2}, \ZZ_T^{0} }
 ds\\&=&\langle  \w p(\cdot,0),\psi\rangle_{ \ZZ_T^{-2}, \ZZ_T^{0} }
 +\int_0^t \langle \w p(\cdot,s),\w  A(\cdot,s)\psi\rangle_{ \ZZ_T^{-2}, \ZZ_T^{0} }
  ds.
\eaaa
Hence
\baaa \langle  \w p'_t(\cdot,t),\psi\rangle_{ \ZZ_T^{-2}, \ZZ_T^{0} } =
\langle  \w p(\cdot,t),\w A(\cdot,t)\psi\rangle_{ \ZZ_T^{-2}, \ZZ_T^{0} }
\quad \hbox{for a.e.}
\quad t.
\eaaa

Let us show that $\w p=\L^* (\w u(\cdot))\rho$. For this, it suffices to show that
\baaa
(\xi,\w p)_{X^0}=
\langle  \rho,\oo V(\cdot,0)\rangle_{ H^{-1}, H^{1} }
\quad \hbox{for any}\quad\xi \in X^0,\quad \oo V=L(\w u(\cdot)) \xi.
\eaaa

For this  $\oo V$, we have that
\baaa
&&\langle  \rho,\oo V(\cdot,0)\rangle_{ H^{-1}, H^{1} }
=\langle  \rho,\oo V(\cdot,0)\rangle_{ H^{-1}, H^{1} }-
\langle \w p(\cdot,T),\oo V(\cdot,T)\rangle_{ \ZZ_T^{-2}, \ZZ_T^{0} }
\\&&=-\int_0^T
\Bigl[(\w p(\cdot,t),\oo V'_t(\cdot,t))_{\ZZ_T^0}+\langle \w p'_t(\cdot,t),\oo V(\cdot,t)\rangle_{ \ZZ_T^{-2}, \ZZ^{0} } \Bigr]dt
\\
&&=-\int_0^T
\Bigl[(-\w A\oo V(\cdot,t)-\xi,\w p(\cdot,t))_{\ZZ_T^0}+ \langle \w p'_t(\cdot,t),\oo V(\cdot,t)\rangle_{ \ZZ_T^{-2}, \ZZ^{0} } \Bigr]dt
\\&&
=(\xi,\w p)_{X^0}+\int_0^T
\Bigl[(\w A\oo V(\cdot,t),\w p(\cdot,t))_{\ZZ_T^0}-\langle \w p'_t(\cdot,t),\oo V(\cdot,t)\rangle_{ \ZZ_T^{-2}, \ZZ^{0} }  \Bigr]dt=(\xi,\w p)_{X^0}.
\eaaa
Hence $\w p=\L^* (\w u(\cdot))\rho$.

Further, we have that
\baaa &&F(a,v_i(\cdot))-F(a,\w u(\cdot))=\E\int_Q \bigl(
p_i(x,t)\varphi (x,v_i(x,t),Z(t),t)- \w p(x,t)\varphi (x,\ww
u(x,t),Z(t),t)\bigr)dx \,dt
\\
&&=\E\int_Q \bigl(  \w p(x,t)\varphi (x,v_i(x,t),Z(t),t)-\w p(x,t)\varphi
(x,  \w u(x,t),Z(t),t))dx\,dt
\\
&&+\E\int_Q  \bigl(    p_i(x,t)\varphi (x,v_i(x,t),Z(t),t)-\w p(x,t)\varphi
(x,v_i(x,t),Z(t),t))dx\,dt \to 0 \ \ \ \hbox{as}\ \ i \to +\infty .
\eaaa
Hence    $F(a,\w
u(\cdot))=J$. This proves the existence of an optimal
control  $\w u(\cdot)\in
U_R$ for problem (\ref{P}) and completes the proof of Theorem \ref{ThEr}.
$\Box$

\par
\subsection{Proof of Theorem  \ref{ThM}.}   Let $a\in\AAA$
be such that $a$ has the probability density function $\rho\in H^0$  such that $\rho(x)>0$ for any $x\in D$.
Let $\w u\in U_R$ be the optimal control that exists by Lemma \ref{ThEr}.
For a given  $\mu\in  U_R$ and $\e\in[0,1]$,  we consider a family   of controls $u_\e=u_{\e,\mu}\in U_R$ such that $
 u_\e =(1-\e  )\w u+\e   \mu (\cdot)$ for all $\e \in [0,1]$.
Since $\Delta_R$ is a convex set,  we have that $u_\e \in U_R$ $(\forall \e ,\mu )$.
\par
We denote
\baaa
&&A(x,t)\defi A(x,\w u(x,t),Z(t),t),\qquad A_\e (x,t)\defi A(x,u_\e (x,t),Z(t),t),\quad
\\&& \oo V_\e (x,t)\defi \oo L(u_\e(\cdot))\varphi(\cdot,u_\e(\cdot),Z(\cdot),\cdot),\quad  V_\e (x,t)\defi L(u_\e(\cdot))\varphi(\cdot,u_\e(\cdot),Z(\cdot),\cdot),
\\&&
\oo{\w V}(x,t)\defi \oo L(\w u(\cdot))\varphi(\cdot,\w u(\cdot),Z(\cdot),\cdot),
\\&&\w\varphi\defi \varphi(\cdot,\w u(\cdot,t),Z(\cdot),\cdot), \quad
\varphi_\e\defi \varphi(\cdot,u_\e(\cdot,t),Z(\cdot),\cdot),\quad
\ww \varphi_\e\defi - \frac{\p\oo V_\e}{\p
t}-\w A\oo V_\e.
\eaaa
Let $\Phi(u(\cdot))\defi F(a, u(\cdot))$.

\begin{proposition}\label{prop3.1} For any
$\e\in[0,1]$,  we have that
\baa &&\Phi  (u_\e (\cdot))-\Phi (\w
u(\cdot))\nonumber
\\&&=\E\int_Q \w p(x,t)[(A_\e (x,t)-\w
A(x,t))\w V(x,t)
     +\varphi             (x,u_\e (x,t),Z(t),t)
-\varphi (x,\w u(x,t),Z(t),t)]dx\,dt\nonumber \\
&& \hphantom{xxxxxxxxx}+\E \zeta_\e,\label{3.13}\eaa
where
\baa
\zeta_\e(x,t)\defi \int_Q \w
p(x,t)(A_\e (x,t)- \w A(x,t))(V_\e (x,t)-\w V(x,t))dx\,dt.
 \label{xi}\eaa
\end{proposition}

{\em Proof of Proposition \ref{prop3.1}}.
It follows from the definitions that
\baaa
 \Phi \bigl(\w u(\cdot)\bigr)=(\w V(\cdot,0),\w p(\cdot,0))_{\ZZ_T^0}=\int_0^T(\w \varphi,\w p)_{\ZZ_T^0}dt,
 \eaaa
 and
 \baaa
 \Phi \bigl(u_\e(\cdot)\bigr)=(V_\e(\cdot,0),\w p(\cdot,0))_{\ZZ_T^0}=
 \int_0^T(\ww \varphi_\e,\w p)_{\ZZ_T^0}dt= \int_0^T(\varphi_\e,\w p)_{\ZZ_T^0}dt
 + \int_0^T((A_\e-\w A)\oo V_\e,\w p)_{\ZZ_T^0}dt,
  \eaaa
since, by the definitions,  \baaa
\varphi_\e=
- \frac{\p \oo V_\e}{\p t}-A_\e \oo V_\e=\ww\varphi_\e-(A_\e-\w A)\oo V_\e.
\eaaa
Using that $\w p \in X^0$, we obtain that
   \baaa
 \Phi \bigl(u_\e(\cdot)\bigr)= \int_0^T(\varphi_\e,\w p)_{\ZZ_T^0}dt
 + \int_0^T((A_\e-\w A)V_\e,\w p)_{\ZZ_T^0}dt.
 \eaaa
 Hence
 \baaa
 \Phi \bigl(u_\e(\cdot )\bigr) -\Phi \bigl( \w u(
\cdot)\bigr)
=\int_0^T
\left\{ (A_\e \w V ,\w p)_{\ZZ_T^0} -(\w A\w V,\w p)_{\ZZ_T^0}
+(\varphi_\e -\w \varphi, \w p)_{\ZZ_T^0} \right\}dt+\E\xi_\e,
\eaaa
since   \baaa
&&
\zeta_\e=\int_0^T\left\{ ((A_\e -\w A )(V_\e-\w V) ,\w p)_{\ZZ_T^0}\right\}dt
\\&&
=\int_0^T
\left\{ (A_\e V_\e ,\w p)_{\ZZ_T^0} -(\w A V_\e ,\w p)_{\ZZ_T^0} -(A_\e \w V,\w p)_{\ZZ_T^0} -(\w A \w V,\w p)_{\ZZ_T^0}
 \right\}dt .
\eaaa

\begin{proposition}\label{lemma3.1} There exists a limit
\baa
\lim_{\e\to 0}\frac{\Phi
(u_\e (\cdot))-\Phi   (\w u(\cdot))}{\e} =\E\int_Q\w
p(x,t)\biggl\{A(x,\mu (x,t),Z(t),t)\w V(x,t)+\varphi (x,\mu (x,t),Z(t),t)
\nonumber\\
 -A(x,\w u(x,t)),t)\w
V(x,t) -\varphi (x,\w u(x,t),Z(t),t)\biggr\}dx\,dt.
\label{3.14}\eaa
\end{proposition}

\par
{\em Proof of Proposition \ref{lemma3.1}}.  By the choice of the a family   of controls,
the first integral in the right hand size of (\ref{3.13})
coincides with the right hand side of (\ref{3.14}) multiplied by $\e$, for any $\e \in (0,1]$.

Let  $W_\e \defi V_\e -\w V$. The lemma will be proved if we show
that
\baaa
J_\e =\e^{-1}\E\int_Q \w p(x,t) \left(A_\e (x,t)-\w A(x,t) \right) W_\e (x,t)dx\,dt \to 0 \
\hbox{as} \ \  \e \to 0.
\eaaa

Let
\baaa
\w b_\e(x,t)=    b(x,u_\e (x,t),Z(t),t)-b(x,\w u(x,t),Z(t),t)=\e \bigl[ b(x,\mu     (x,t),t)-b(x,\w u(x,t),Z(t),t)\bigr],
\\
\w f_\e (x,t)=  f(x,u_\e (x,t),Z(t),t)-f(x,\w u(x,t),Z(t),t)=\e \bigl[f(x,\mu    (x,t)-f(x,\w
u(x,t),Z(t),t)\bigr],
\\
\w\varphi_\e (x,t)=  \varphi(x,u_\e (x,t),Z(t),t)-\varphi(x,\w u(x,t),Z(t),t)=\e \bigl[\varphi(x,\mu    (x,t)-\varphi(x,\w
u(x,t),Z(t),t)\bigr].
\Lindex{\\
\w \lambda_\e (x,t)= \lambda(x,u_\e (x,t),t)-\lambda(x,\w u(x,t),t)=\e
 \bigl[\lambda(x,\mu (x,t),t)
-\lambda(x,\w u(x,t),t)\bigr]}
\eaaa
 The second equalities in the above formulae follow from the choice of the a family   of controls.

Let \baaa
 \xi_\e =-\w
\varphi_\e (x,t) +\sum_{i,j=1}^n \w b^{(\e)}_{ij}(x,t) \frac{\p^2
\w V}{\p x_i \p x_j}(x,t)+ \frac{\p      \w V}{\p
x}(x,t) \w f_\e (x,t)\Lindex{-\w V(x,t)\w \lambda_\e (x,t)}. \eaaa

By the choice of $u_\e$, it follows that
\baaa
\E|\xi_\e|\le \const \e \|W_\e\|_{\W^2}.
\eaaa

       \par Further,  $W_\e $ is such that, in  $Q$,
\baaa
&&\frac{\p W_\e }{\p  t}  +A_\e W_\e  =\xi_\e ,\\
&&W_\e |_{x\in\p D}=0,\quad W_\e |_{t=T}=0,
\eaaa
We have that
$
| J_\e | \le C\|  W_\e \|_{Y^2}\| \w p\|_{L_2(Q)} \le C\|\xi_\e \|_{X^0}\| \rho \|_{H^{-1}  }
$
for a constant  $C>0$.  Since  $b,f,\lambda $ are bounded functions, we have that
\baaa
 |\xi_\e (x,t)|\le\e
C_1\left(\sum_{i,j=1}^{n}| \frac{\p^2 \w V}{\p x_i \p  x_j}|  +|
\frac{\p  \w V}{\p x}(x,t)| +\Lindex{| \w V(x,t)| +}\sup_{v\in \Delta}|
\varphi (x,v,t)| \right), \eaaa where $C_1 >0$.

From the choice  of the controls
$u_\e (\cdot)$, we obtain that  $\| \xi_\e\|_{X^0}\le \e \hbox{  const} \| \w
V\|_{Y^2}.$  This completes the proof of   Proposition \ref{lemma3.1}.
\par We are now in the
position to prove Theorem \ref{ThM}(i).
\par
Let  $\Xi $ be the set of     $(x,t)$ that are Lebesgue points of  $$\w p(x,t)[A(x,\w u(x ,t),Z(t),t)\w V(x,t)+\varphi (x,\w u(x,t),Z(t),t)].
$$
Clearly,  $\hbox{ mes} \{(\Theta \times
[0,T]) \backslash \Xi  \}=0$.   By the continuity
$b(x,v,t)$, $f(x,v,t)$, \Lindex{$\lambda (x,v,t)$,} and  $\varphi (x,v,t)$, by  Luzin Theorem
 from Shilov and Gurevich (2012),
 \index{\cite{ShilovG}} p.87,  we obtain that for all
$v\in \Delta $  the vectors   $(x ,t)\in \Xi $ are Lebesgue points of
$\w p(x,t)(A(x ,v,t)\w V(x,t)+\varphi
(x ,v,t))$.
\par
Further, let us consider
 $\mu (\cdot)\in U_R$  such that $\mu
(x ,t)=\w u(x ,t)$ for $(x ,t)\notin B$,  $\mu
(x ,t)=v$ for $(x ,t)\in B$, where $v\in \Delta $
represent Dirac measure,
$B\subset \Theta  \times (0,T)$ are arbitrary domains such that they form a
 Vitali system of sets properly shrinking in the sense of definition from Shilov and Gurevich (2012)
  \index{ShilovG}
 to each point
$(x,t)\in \Xi $. Hence
\baaa
A(x,\w u(x,t),t)\w V(x,t)+\varphi (x,\w u(x,t),Z(t),t) \le  A(x,v,t)\w
V(x,t)+\varphi (x,v,Z(t),t) \quad \hbox{a.e.}\quad \forall v\in \Delta_R.
\eaaa
Then the statement of  the proof of Theorem  \ref{ThM}(i)  for $\w u\in U_R$
 follows from Proposition \ref{lemma3.1}.  To complete
 the proof of Theorem  \ref{ThM}(i), we need to show that   there exists $\w u\in U_0$
 with the required properties.

 \begin{proposition}\label{prop3.5}
There exists $\ww u(\cdot)\in U_0$ such that, for a.e.   $x,t$
\baaa
A(x,\ww u(x,t),Z(t),t)\w V(x,t)+\varphi (x,\ww u(x,t),Z(t),t) \le  A(x,v,Z(t),t)\w
V(x,t)+\varphi (x,v,Z(t),t) \quad \forall v\in \Delta.
\eaaa
and
\baa
\w V=L(\w u(\cdot))\varphi(\cdot,\w u(\cdot),Z(\cdot),\cdot)  =L(\ww
u(\cdot))\varphi (\cdot,\ww u(\cdot),Z(\cdot),\cdot),
\label{LL}\eaa
where  $\w u\in U_R$ is an optimal control for described in the proof above.
 \end{proposition}
  \par

{\em Proof of  Proposition \ref{prop3.5}}.
 Let  $R\defi D\times \Delta
\times [0,T]$, and let  $z(x,v,t)\defi A(x,v,t)\w V(x,t)+\varphi
(x,v,t)$. Let $R'\subset R$ be such that  $\hbox{ mes}
\{R\backslash R'\}=0$, and $R'=\cup_{k=1}^{+\infty} R_k$, where  $R_k=R_k(\o)$
are random $\F^Z_T$-measurable  compact sets defined a.s. such that the function  $z(x,v,t)$ is continuous on
 $R_k$.  (These  $R_k=R_k(\o)$  exist by the Luzin Theorem   from
 Shilov and Gurevich (2012), \index{\cite{ShilovG},} p.87.
 Let   \baaa
 S_k\defi \{(x,\w v,t)\in R_k:
z(x,\w v,t)=\inf_{v\in \Delta}z(x,v,t)\},\quad S=\cup_{k=1}^{+\infty} S_k.
\eaaa
The set  $S$ is $\sigma $-compact.
 By Lemma  B from Fleming an Rishel (1975), \index{ \cite{FlemRish},} p.277,
there exists a desired function     $\ww
u(\cdot)\in U_0$ such that  $(x,\ww u(x,t),t)\in S$ a.s. for a.e.
$x,t$.
 In particular, this means that
\baaa
A(x,\ww u(x,t),t)\w V(x,t)+\varphi (x,\ww u(x,t),Z(t),t) =A(x,\w
u(x,t),t)\w V(x,t)+\varphi (x,\w u(x,t),Z(t),t)\quad \hbox{a.e.}. \eaaa
hence (\ref{LL}) holds.
This completes the proof of Proposition \ref{prop3.5}. $\Box$

 The proof of Theorem  \ref{ThM}(i) follows from Proposition \ref{prop3.5}.

 Let us prove Theorem \ref{ThM}(ii). Let $\w a\in\AAA$, and let
 $\w u\in U_0$ be an  optimal control for corresponding problem (\ref{P})
  that exists according  to Theorem \ref{ThM}(i). Let
  $\w V=L(\w u(\cdot))\varphi(\cdot,\w u(\cdot),Z(\cdot),\cdot)$.

  Further,   let  $a\in\AAA$ be such that its probability density
  $\rho\in H^0$. Let
   $ u\in U_R$  be any.
 It follows from the definitions   that
  \baaa
  &&A(x,\w u(x,t),Z(t),t)\w V(x,t)+\varphi(x,\w u(x,t),Z(t),t)
 \\ &&=A(x, u(x,t),Z(t),t)\w V(x,t)+\varphi(x,u(x,t),Z(t),t)
 +\psi(x,t),
  \eaaa
  where
  \baaa
  \psi(x,t)\defi A(x,\w u(x,t),Z(t),t)\w V(x,t)+\varphi(x,\w u(x,t),Z(t),t)\\- A(x, u(x,t),Z(t),t)\w V(x,t)-\varphi(x,\w u(x,t),Z(t),t).
  \eaaa
Hence   $\w V=L(u(\cdot))[\varphi(\cdot,u(\cdot),\cdot)+\psi]$. By the definition
of the operator $L*(u(\cdot))$, we have that
\baaa
  (\w V(\cdot,0),\rho)_{\ZZ_T^0}=\E\int_0^T(\varphi(\cdot,\w u(\cdot),Z(\cdot),t),p_u(\cdot,t))_{\ZZ_T^0}dt
  +\oo R,
  \eaaa
  where $p_u\defi L(u(\cdot))^*\rho$,
  \baaa
  \oo R \defi \E\int_0^T(\psi(\cdot,t),p_u(\cdot,t))_{\ZZ_T^0}dt.
  \eaaa
On the other hand, it follows from the definitions that
\baaa
(\w V(\cdot,0),\rho)_{\ZZ_T^0}=F(a,\w u(\cdot)),\qquad \E\int_0^T(\varphi(\cdot,\w u(\cdot),t),p_u(\cdot,t))_{\ZZ_T^0}dt=F(a,u(\cdot)).
\eaaa

By the choice of  $p_u$, $\w u$, and $\w V$,
  \baaa
 p_u(x,t,\o)\ge 0,\quad  \psi(x,t,\o)\le 0\quad \hbox{a.e.}
 \label{psi+} \eaaa
Hence  $\oo  R\le 0$.  This proves  the statement  of  Theorem \ref{ThM}(ii) and  completes the proof of  Theorem \ref{ThM}. $\Box$
\section{Discussion and further research}
Similarly to the first order SPDEs introduced in Bender and Dokuchaev (2016a,b),
equation  (\ref{HJB}) is an analog of the HJB equation. In our case, it has some special features. For example, solutions of (\ref{HJB}) are not necessarily   continuous in $t$ since the filtration $\F^Z_t$ can be discontinuous for the case where $Z(t)$
is discontinuous.

The proofs in the present paper are different from the proofs from
 Bender and Dokuchaev (2016a,b) and from the proofs from Dokuchaev (2017); the present proof
 since it uses the regularity properties of non-degenerate parabolic equations
and does use the time discretisation implemented in the cited papers. Therefore,
the proofs in the present  are  more straightforward.
However, it is unlikely that this proof
can be extended on the special  control problems  considered in Bender and Dokuchaev (2016a,b) and  Dokuchaev (2017);
these problems can be regarded as degenerate problems in domains with boundaries.
The regularity of corresponding value functions in these papers  was not covered by the existing literature;  it was analyzed  directly using the time discretisation.

Alternatively to solution of equation  (\ref{HJB}), the value function $\w V$ can be estimated via Monte-Carlo method combined with the dual pathwise optimization method, similarly to Section  7 in Bender and Dokuchaev (2016a) or Section 4 in Dokuchaev (2017).
In this case,  equation  (\ref{HJB}) can be useful of calculation of the optimal control  as the process where the minimum in Proposition \ref{prop3.5} is achieved for optimal $\w V$.

The present paper considers only the case where the value function has $L_2$-integrable
second order derivatives in $x$.
It could be interesting to extend the results of this this paper on more general class diffusion coefficients. This may require to consider viscosity  solutions of parabolic equations
with measure-valued derivatives.


\begin{thebibliography}{100}
\bibitem{Alos}
E.  Al\'os, J. A. Le\'on, D. Nualart. (1999).
 Stochastic heat equation with random coefficients
 {\it
Probability Theory and Related Fields}, {\bf 115} (1), 41--94.

\bibitem{Ba}
V.  Bally, I. Gyongy, E. Pardoux. (1994).  White noise driven
parabolic SPDEs with measurable drift, {\it Journal of Functional
Analysis} {\bf 120}, 484--510.

\bibitem{BD1} C. Bender and  N. Dokuchaev. (2016a). A
first-order BSPDE for swing option pricing.   {\em
Mathematical Finance} {\bf 26} (3), 461--491.
\bibitem{BD2} C. Bender and  N. Dokuchaev. (2016b). A
first-order BSPDE for swing option pricing: Classical solutions.    {\em
Mathematical Finance}, online published in 2015, in press.

\bibitem{Cordes}
H.O. Cordes. (1956). Uber die erste
Randwertaufgabe
bei quasilinerian
Differentialgleichun-genzweiter Ordnung in mehr als zwei Variablen, Math. Ann. 131, 278-312.



\bibitem{CK} T.  Caraballo,  P. E. Kloeden,
B. Schmalfuss. (2004).  Exponentially stable stationary solutions for
stochastic evolution equations and their perturbation, {\em Appl. Math.
Optim.} {\bf  50}, 183--207.
\bibitem{CMG}
A. Chojnowska-Michalik, B. Goldys. (1995). {Existence,
uniqueness and invariant measures for stochastic semilinear
equations in Hilbert spaces},  {\it Probability Theory and Related
Fields},  {\bf 102}(3), 331--356.


\bibitem{CL}
M. Crandall and P. L. Lions. (1983a). Viscosity solutions of Hamilton-Jacobi equations,
Trans. Amer. Math. Soc. 277, 1-42

\bibitem{DaPT}
G. Da Prato,   L. Tubaro. (1996). { Fully nonlinear stochastic
partial differential equations}, {\it SIAM Journal on Mathematical
Analysis} {\bf 27}, No. 1, 40--55.
\bibitem{D97} N.G. Dokuchaev. (1997). Cordes conditions and some alternatives
for parabolic equations and discontinuous diffusion. {\it
Differential equations } {\bf 33} , N 4, 433-442.

\bibitem{D05}  Dokuchaev, N. (2005). Parabolic
Ito equations and second fundamental inequality.  {\it Stochastics}
{\bf 77}, iss. 4., pp. 349-370.

\index{\bibitem{D08b}
N. Dokuchaev (2008). Parabolic Ito equations with mixed in time
conditions,
{\it Stochastic Analysis and Applications},  {\bf 26} (3), 562--576. 
}

\bibitem{D18}
N. Dokuchaev. (2018).
On degenerate backward SPDEs in bounded domains under non-local conditions.
{\em Stochastics} 90 (8), 1170-1189.
   \bibitem{DT}
K. Du,  S. Tang. (2012). Strong solution of backward stochastic
partial differential equations in $C^2$ domains. {\em  Probability
Theory and Related Fields},  {\bf 154},  255--285.
 \bibitem{DTZ}
K. Du, S. Tang, Q. Zhang. (2013). $W^{m,p}$-solution ($p\ge 2$) of linear degenerate backward stochastic partial differential equations in the whole space.
{\em Journal of Differential Equations} {\bf 254} (7), 2877--2904.
 \bibitem{DZ}
K. Du, Q. Zhang. (2013).
Semi-linear degenerate backward stochastic partial differential equations and associated forward–backward stochastic differential equations
{\em Stochastic Processes and their Applications} {\bf 123} (5), 1616--1637.
\bibitem{Duan}
J. Duan,  K. Lu, B. Schmalfuss. (2003). Invariant manifolds for
stochastic partial differential equations, {\em Ann. Probab.},  {\bf
31},
 2109--2135.

\bibitem{FZ}
C. Feng, H. Zhao. (2012). Random periodic solutions of SPDEs via
integral equations and Wiener-Sobolev  compact embedding. {\em
Journal of Functional Analysis}  {\bf 262}, 4377--4422.

\bibitem{FlemRish}
W.H. Fleming,  R.W. Rishel. (1975).  Deterministic and Stochastic Optimal Control. New York-Heidelberg-Berlin. Springer-Verlag.
\bibitem{FlemSoner}
W. Fleming and H. M. Soner. (1993). Controlled Markov Processes and Viscosity Solutions,
Springer, Berlin-Heidelberg-New York.

\bibitem{Gv2}
I.I. Gihman and A.V. Skorohod. (1975). {\it The Theory of
Stochastic Processes.} Vol. 2. Springer-Verlag, New York.

 \bibitem{GM}
I. I. Gikhman, T.M. Mestechkina. (1983). The Cauchy problem for stochastic first-order partial differ- ential equations. Theory of Random Processes 11, 25--28.
\bibitem{A1}
I. Gy\"ongy. (1998). Existence and uniqueness results for semilinear
stochastic partial differential equations. {\it Stochastic Processes
and their Applications},  {\bf 73} (2),  271--299.


\bibitem{HuPeng}
Y. Hu, S. Peng. (1995).
Solution of forward-backward stochastic differential equations. {\em
Probability Theory and Related Fields} {\bf103},  273--283.

 \bibitem{Hu} Y. Hu, J. Ma, J. Yong. (2002).   On semi-linear degenerate backward stochastic partial differential equations,
{\em  Probab. Theory Related Fields} {\bf 123} (3), 381--411.

\bibitem{Kalita}
E.A. Kalita. (1989). Regularity of solutions of
Cordes-type elliptic systems of any order. {\it Doklady Acad. Sci.
Ukr. SSR, A}, {\bf  5}, 12-15.
\bibitem{Ko}
 A.I. Koshelev. (1982). On exact conditions of regularity of
solutions of for elliptic systems and  Liouville theorem. {\it
Dokl. Akad. Nauk. SSSR}, 265, Iss.6, 1309-1311.

\bibitem{KR}
N.V. Krylov.   (1980).  {\em Controlled diffusion processes}.
New York, USA: Springer,.
\bibitem{KR2}
N.V. Krylov. (1987). {\em Nonlinear Elliptic and Parabolic Equations of the Second Order}.
Springer, Netherlands.
\bibitem{Ku}
H. Kunita. (1990) Stochastic Flows and Stochastic Differential Equations. Cambridge University Press, Cambridge.

\bibitem{la} O.A.
Ladyzhenskaia.  (1985).
{\it The Boundary Value Problems of Mathematical Physics}. New York:
Springer-Verlag.
 \bibitem{Landis}
E.M. Landis. (1998).  {\it Second Order equations of elliptic and
parabolic Type}, vol. 171, Amer. Math. Soc., Providence, R.I.,
English transl. in Translations of Math. Monographs.


\bibitem{V4}
P.L. Lions. (1982). Generalized solutions of Hamilton-Jacobi equations. Research
Notes in Mathematics, Vol. 69, Pitman Advanced Publishing Program,
Boston, 317 pp. 
\bibitem{liu}
Y.  Liu, H.Z. Zhao. (2009). Representation of pathwise stationary
solutions of stochastic Burgers equations,  {\em Stochastics and
Dynamics}  {\bf  9} (4), 613--634.
\bibitem{Ma99}
 J. Ma, J. Yong, J. (1999). On linear, degenerate backward
stochastic partial differential equations.
 {\em Probability Theory and Related Fields}  {\bf 113} (2) (1999), 135--170.
 \bibitem{Mas}
B. Maslowski. (1995). { Stability of semilinear equations with
boundary and pointwise noise}, {\it Annali della Scuola Normale
Superiore di Pisa - Classe di Scienze} (4),  {\bf 22},  No. 1, 55--93.
\bibitem{Mat}
J. Mattingly. (1999). Ergodicity of 2D Navier-Stokes equations with
random forcing and large viscosity, {\em Comm. Math. Phys.},  206 (2),
273--288.
\bibitem{Moh}
S.-E. A.  Mohammed, T.  Zhang,  H. Z. Zhao. (2008). The stable manifold
theorem for semilinear stochastic evolution equations and stochastic
partial differential equations, {\em Mem. Amer. Math. Soc.} {\bf 196}
(917),  1-105.



\bibitem{Par}
E. Pardoux. (1993).  Stochastic partial differential equations, a
review, {\em Bull. Sci. Math.} {\bf 117} (1), 29--47.

\bibitem{Peng}
S. Peng. (1992). Stochastic Hamilton-Jacobi-Bellman equations.  {\em  SIAM J. Control Optim.} 30, 284--304.

\bibitem{Roz}
B. L. Rozovskii. (1990).  {\it Stochastic Evolution Systems; Linear
Theory and Applications to Non-Linear Filtering,}  Kluwer Academic
Publishers,  Dordrecht-Boston-London.

\bibitem{ShilovG}
G.E. Shilov and B. L. Gurevich. (2012).  Integral, Measure, and Derivative.
Dover Books on Mathematics. New York.
\bibitem{Tal}
G. Talenti. (1965). Sopra una classe di equazioni
elilitticche a coefficienti misurabili. {\it Ann. Math. Pure.
Appl.} {\bf 69}, 285-304.

\bibitem{Varga}
	J. Varga. (1972). Optimal Control for Differential and functional Equations.
Academic Press, New York.


  \bibitem{Wa}
J. B. Walsh. (1986). An introduction to stochastic partial
differential equations, Lecture Notes in Mathematics, {\bf 1180},
Springer, New York.

\bibitem{Y}
 K. Yosida, Functional Analysis (Springer, Berlin Heilderberg New York, 1995.


\bibitem{Zh}
X. Y. Zhou. (1992). { A duality analysis on stochastic partial
differential equations}, {\it Journal of Functional Analysis}  {\bf
103} (2), 275--293.

\end{thebibliography}
\end{document}